\documentclass[11pt]{article}
\usepackage{a4}
\usepackage{latexsym}
\usepackage{amsfonts}
\usepackage{amssymb}
\usepackage{amscd}
\usepackage{amsmath}
\usepackage{epsfig}
\begin{document}
\title{\textbf{On the poles of topological zeta functions}}
\author{Ann Lemahieu, Dirk Segers and Willem Veys\footnote{Ann Lemahieu, Dirk Segers, Willem Veys, K.U.Leuven, Departement Wiskunde,
Celestijnenlaan 200B, B-3001 Leuven, Belgium, email:
ann.lemahieu@wis.kuleuven.ac.be, dirk.segers@wis.kuleuven.ac.be,
wim.veys@wis.kuleuven.ac.be} \date{}}

\maketitle {\footnotesize \emph{\textbf{Abstract.---} We study the
topological zeta function $Z_{top,f}(s)$ associated to a
polynomial $f$ with complex coefficients. This is a rational
function in one variable and we want to determine the numbers that
can occur as a pole of some topological zeta function; by
definition these poles are negative rational numbers. We deal with
this question in any dimension. Denote $\mathcal{P}_n := \{ s_0
\mid \exists f \in \mathbb{C}[x_1,\ldots, x_n] \, : \,
Z_{top,f}(s) \textsl{\mbox{ has a pole in }} s_0 \}$. We show that
$\{-(n-1)/2-1/i \mid i \in \mathbb{Z}_{>1}\}$ is a subset of
$\mathcal{P}_n$; for $n=2$ and $n=3$, the last two authors proved
before that these are exactly the poles less then $-(n-1)/2$. As
main result we prove that} each \emph{rational number in the
interval $[-(n-1)/2,0)$ is contained in $\mathcal{P}_n$.}}
 \\ ${}$
\begin{center}
\textsc{1. Introduction}
\end{center} ${}$\\
\indent Denef and Loeser created in $1992$ a new zeta function,
which they called the topological zeta function because of the
topological Euler--Poincar\'e characteristic turning up in it.
Roughly said, the topological zeta function $Z_{top,f}$ associated
to a polynomial $f$ is a function containing information we can
pick out of each chosen embedded resolution of $f^{-1} \{ 0 \}
\subset \Bbb A^n$. They introduced it in \cite{DenefLoeser1} in
the following way.
\\ \indent Let $f$ be a polynomial in $n$ variables over $\mathbb{C}$ and let $h : X \rightarrow \Bbb A^n$ be an
embedded resolution of $f^{-1} \{ 0 \}$. To define $Z_{top,f}$ we
need some data related to the embedded resolution $(X,h)$. Let
$E_i, i \in S$, be the irreducible components of $h^{-1}(f^{-1}\{
0 \})$, then denote by $N_i$ and $\nu_i - 1$ the multiplicities of
$E_i$ in the divisor on $X$ of $f \circ h$ and $h^\ast (dx_1
\wedge \ldots \wedge dx_n)$, respectively. The couples $(N_i,
\nu_i), i \in S$, are called the numerical data of the resolution
$(X,h)$. For $I \subset S$ we denote also $E_I := \cap_{i \in I}
E_i$ and $E^\circ_I := E_I \setminus (\cap_{j \notin I} E_j)$.
Further we write $\chi(\cdot)$ for the topological
Euler--Poincar\'e characteristic.
\\ \\
\noindent \textbf{Definition.---} The local \emph{topological zeta
function associated to $f$} is the rational function in one complex
variable\[Z_{top,f} (s)  := \sum_{I \subset S} \chi (E^\circ_I \cap
h^{-1}\{ 0 \}) \prod_{i \in I} \frac{1}{N_i s+ \nu_i}.\] There is a
global version replacing $E^\circ_I \cap h^{-1}\{ 0 \}$ by
$E^\circ_I$. When we do not specify, we mean the local one.
\\ \indent Denef and Loeser proved that every embedded resolution
gives rise to the same function, so the topological zeta function
is a well-defined singularity invariant (see \cite{DenefLoeser1}).
Once the motivic Igusa zeta function was introduced, they proved
this result alternatively in \cite{DenefLoeser2} by showing that this more general zeta
function specialises to the topological one.
\\ \indent In particular the poles of the topological zeta function of
$f$ are interesting numerical invariants. Various conjectures
relate them to the eigenvalues of the local monodromy of $f$, see
for example \cite{DenefLoeser1}. The poles are part of the set
$\{-\nu_i/N_i \mid i \in S\}$; therefore the
$-\nu_i/N_i$ are called the candidate poles. Notice that
the poles are negative rational numbers.
\\ \indent A related numerical invariant of $f$ at $0 \in \mathbb{C}^n$
is its \emph{log canonical threshold} $c_0(f)$ which is by
definition \[\sup \{ c \in \mathbb{Q} \mid \mbox{ the pair }
(\mathbb{C}^n , c \; \mbox{div} \; f) \mbox{ is log canonical in a
neighbourhood of } 0 \}.\] It is described in terms of the
embedded resolution as $c_0(f)=\min\{ \nu_i/N_i \mid 0 \in h(E_i),
i \in S \}$ (see \cite[Proposition 8.5]{Kollar2}). It was studied
in various papers of Alexeev, Cheltsov, Ein, de Fernex, Koll\'ar,
Kuwata, M$^{c}$Kernan, Musta\c t\u a, Park, Prokhorov, Reid,
Shokurov and others. Especially the sets $\mathcal{T}_n := \{
c_0(f) \mid f \in \mathbb{C}[x_1,\ldots, x_n] \},$ with $n \in
\mathbb{Z}_{>0}$, show up in interesting conjectures, see
\cite{Alexeev}, [Ko], [Ku],
\cite{M$^c$KernanProkhorov}, [Pr] and \cite{Shokurov}. \\
\indent In the context of the topological zeta function, it is
natural to study similarly the set
\[ \mathcal{P}_n := \{ s_0 \mid \exists f \in \mathbb{C}[x_1,\ldots, x_n] \,
: \, Z_{top,f}(s) \textsl{\mbox{ has a pole in }} s_0 \}. \] The
case $n=1$ is trivial: $\mathcal{P}_1= \{-1/i \mid i \in
\mathbb{Z}_{>0}\}$. \\ \indent From now on we assume that $n \geq
2$. A more or less obvious lower bound for $\mathcal{P}_n$ is
$-(n-1)$, see \cite[Section 2.4]{Segersthesis}. In
\cite{SegersVeys}, the second and the third author studied the
`smallest poles' for $n=2$ and $n=3$. They showed that
$\mathcal{P}_2 \cap (-\infty,-\frac{1}{2})  =
\{-\frac{1}{2}-\frac{1}{i} \mid i \in \mathbb{Z}_{>1} \}$
and that 
$\mathcal{P}_3 \cap (-\infty,-1)  = \{-1-\frac{1}{i} \mid i \in
\mathbb{Z}_{>1} \}$.
They expected that this could be generalised to
\begin{eqnarray*} \mathcal{P}_n \cap (-\infty,-\frac{n-1}{2}) =
\{-\frac{n-1}{2}-\frac{1}{i}\mid i \in \mathbb{Z}_{>1} \},\quad
\mbox{for all } n \in \mathbb{Z}_{>1}. \end{eqnarray*} In
particular, they predicted that the lower bound $-(n-1)$ could be
sharpened to $-n/2$. This better bound was recently proven by the
second author in \cite{Segers}. In this article we verify for all $n
\geq 4$ that $\{-(n-1)/2-1/i \mid i \in \mathbb{Z}_{>1}\}  \subset
\mathcal{P}_n$, and as main result we show that \emph{any} rational
number in the remaining interval
$[-(n-1)/2,0)$ is a pole of some topological zeta function. \\ \\
\textbf{Theorem.---} \emph{For $n \geq 2$ we have
$[-(n-1)/2,0) \cap \mathbb{Q} \subset \mathcal{P}_n$.}
\\ \\ \indent With the Thom-Sebastiani principle \cite{DenefLoeser3}, $x_1^i+x_2^2+\cdots+x_n^2$
is the obvious candidate to have $-(n-1)/2-1/i$ as a pole of its
associated topological zeta function. It is not clear a priori that
this will be true for all $n$ and $i$. We check this in section $2$.
For the theorem, however, the key is to find a suitable family of
polynomials.
\\ \indent We will put the useful information of the resolution into a
diagram, which is called the dual intersection graph. It is
obtained as follows. One associates a vertex to each exceptional
component in the embedded resolution (represented by a dot) and to
each component of the strict transform of $f^{-1}\{0\}$
(represented by a circle). One also associates to each
intersection an edge, connecting the corres\-ponding vertices. The
fact that $E_i$ has
numerical data $(N_i,\nu_i)$ is denoted by $E_i(N_i,\nu_i)$. \\
\indent When the strict transform of $f^{-1}\{0\}$ is irreducible,
we will denote it
by $E_0$. 
Let $E_i$ be an exceptional
variety and let $E_j$, $j \in J$, be the components that intersect
$E_i$ in $X$. We set $\alpha_j:=\nu_j-(\nu_i/N_i)N_j$ for $j \in J$;
these numbers appear in the calculation of the residue of
$Z_{top,f}$ in $-\nu_i/N_i$.
\\ ${}$
\begin{center}
\textsc{2. The set $\{-(n-1)/2-1/i \mid i \in \mathbb{Z}_{>1}\}$ is a subset of $\mathcal{P}_n$}
\end{center} ${}$ \\
\noindent \emph{\textbf{Embedded resolution for
\boldmath{$x_1^i+x_2^2+\cdots+x_n^2=0$}, \boldmath{$n \geq 4$},
with \boldmath{$i$} even}} \\ \indent After blowing up $i/2$ times
in the origin, we get an embedded resolution for $f$. We present
the dual intersection graph for $i \not=2$.
\\ \\
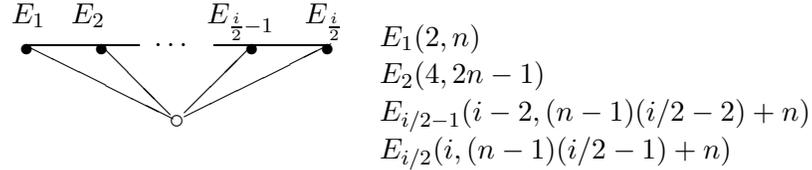
\begin{picture}(150,20)(-5,2)
\put(30,15){\line(1,0){15}}
\put(30,14.5){\circle*{1.5}}
\put(40,14.5){\circle*{1.5}}\put(50,5){\circle{1.5}}
\put(28,18){$E_{1}$} \put(36,18){$E_{2}$}
\put(47,15){$\ldots$}  \put(30,15){\line(2,-1){19.5}}
\put(40,15){\line(1,-1){9.15}} \put(51,6){\line(1,1){9}}
\put(51,5.5){\line(2,1){19.5}} \put(60,14.5){\circle*{1.5}}
\put(70,14.5){\circle*{1.5}}
\put(55,15){\line(1,0){15}}
\put(54,18){$E_{\frac{i}{2}-1}$}
\put(67,18){$E_{\frac{i}{2}}$}
\put(77,15){$E_1(2,n)$} \put(77,10){$E_2(4,2n-1)$}
\put(77,05){$E_{i/2-1}(i-2,(n-1)(i/2-2)+n)$}
\put(77,0){$E_{i/2}(i,(n-1)(i/2-1)+n)$}
\end{picture}  \\
\\The exceptional variety $E_{i/2}$ gives the candidate pole $-(n-1)/2-1/i$ in which we are interested. If $i \not=2$,
its residue is
\begin{eqnarray*}
& \frac{1}{N_{\frac{i}{2}}} & \left( \chi(E^\circ_{I_1}) +
\chi(E^\circ_{I_2})\frac{1}{\alpha_{\frac{i}{2}-1}}+
\chi(E^\circ_{I_3})\frac{1}{\alpha_{0}}+
\chi(E^\circ_{I_4})\frac{1}{\alpha_{0}\alpha_{\frac{i}{2}-1}}
\right),
\end{eqnarray*}
 where $I_1=\{\frac{i}{2}\}, \quad I_2=\{\frac{i}{2},\frac{i}{2}-1\}, \quad I_3=\{\frac{i}{2},0\},
\quad I_4= \{\frac{i}{2},\frac{i}{2}-1,0\}.$ The Euler--Poincar\'e
characteristics $\chi(E^\circ_{I_j})$, $1 \leq j \leq 4$, are put
in Table $1$. These are easily computed since $E_{i/2} \cong
\mathbb{P}^{n-1}$, and $E_{i/2-1}$ and $E_0$ intersect $E_{i/2}$
in a hyperplane and a smooth quadric, respectively.
\begin{center}
\begin{tabular}[t]{|c|c|c|} \hline $\chi(E^\circ_{I_j})$ & n
odd & n even
\\ \hline & & \\
$j=1$ & $1$ & $-1$  \\
$j=2$ & $0$ & $1$\\
$j=3$ & $0$ & $2$ \\
$j=4$ & $n-1$ & $n-2$ \\
& & \\
\hline
\end{tabular}
\end{center}
\begin{center}
\emph{Table 1}
\end{center}
Using that $\alpha_0=(3-n)/2-1/i$ and $\alpha_{i/2-1}=2/i$, some
easy calculations yield that the residue is non-zero, for all $n \in
\mathbb{N}$, $n \geq 4$.
\\ \indent When $i=2$, we blow up just once in the origin to get an embedded resolution. By using
$\alpha_0=\frac{2-n}{2}, \chi(E^\circ_{I_1})=0 (n \mbox{ even}),
\chi(E^\circ_{I_1})=1 (n \mbox{ odd}),$ we conclude that also here
the residue is non-zero.
\\ \\ \noindent \emph{\textbf{Embedded resolution for \boldmath{$x_1^i+x_2^2+\cdots+x_n^2=0$}, \boldmath{$n \geq 4$}, with \boldmath{$i$} odd}} \\
\indent After blowing up $(i+1)/2$ times in the origin, followed by
blowing up once more in $D:=E_{(i+1)/2} \cap E_{(i-1)/2} \cong
\mathbb{P}^{n-2}$, we get an embedded resolution with the following
dual intersection graph.
\\ \\ \\
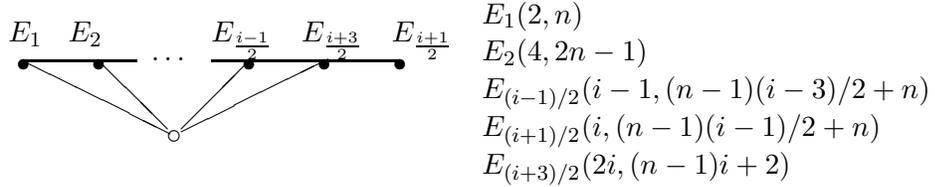
\begin{picture}(150,20)(-5,-5)
\put(15,15){\line(1,0){15}} \put(15,14.5){\circle*{1.5}}
\put(25,14.5){\circle*{1.5}} 
\put(13,17.5){$E_{1}$} \put(21,17.5){$E_{2}$} \put(32,15){$\ldots$}
\put(40,15){\line(1,0){25}} \put(45,14.5){\circle*{1.5}} \put(55,14.5){\circle*{1.5}}
\put(40,17.5){$E_{\frac{i-1}{2}}$} \put(52,17.5){$E_{\frac{i+3}{2}}$}
\put(64,17.5){$E_{\frac{i+1}{2}}$} \put(65,14.5){\circle*{1.5}}
 \put(15,15){\line(2,-1){19.5}}
\put(25,15){\line(1,-1){9.15}} \put(36,6){\line(1,1){9}}
\put(36,5.5){\line(2,1){19.5}}
\put(35,5){\circle{1.5}}
\put(76,20){$E_1(2,n)$} \put(76,15){$E_2(4,2n-1)$}
\put(76,10){$E_{(i-1)/2}(i-1,(n-1)(i-3)/2+n)$}
\put(76,5){$E_{(i+1)/2}(i,(n-1)(i-1)/2+n)$}
\put(76,0){$E_{(i+3)/2}(2i,(n-1)i+2)$}
\end{picture}  \\
The last exceptional variety has $-(n-1)/2-1/i$ as candidate pole.
The rele\-vant subsets in the computation of the residue are
$I_1=\{\frac{i+3}{2}\}, \quad I_2=\{\frac{i+3}{2},0\}, \quad
I_3=\{\frac{i+3}{2},\frac{i+1}{2}\}, I_4=
\{\frac{i+3}{2},\frac{i-1}{2}\},
I_5=\{\frac{i+3}{2},\frac{i-1}{2},0\}.$ Here $E_{(i+3)/2}$ is a
$\mathbb{P}^1$-bundle over $D$. For $j=2,3,4$ we have that
$E_{I_j} \cong D$ and $E_{I_5}$ is a smooth quadric. With the
Euler--Poincar\'e characteristics of Table $2$ and
$\alpha_0=(3-n)/2-1/i$, $\alpha_{(i-1)/2}=1/i$ and
$\alpha_{(i+1)/2}=(n-1)/2$, we find that the residue is non-zero,
for all $n \geq 4$.
\begin{center}
\begin{tabular}[t]{|c|c|c|} \hline $\chi(E^\circ_{I_j})$ & n
odd & n even
\\ \hline & & \\
$j=1$ & $0$ & $-1$  \\
$j=2$ & $0$ & $1$\\
$j=3$ & $n-1$ & $n-1$ \\
$j=4$ & $0$ & $1$ \\
$j=5$ & $n-1$ & $n-2$ \\
& & \\
\hline \end{tabular}
\end{center}
\begin{center}
\emph{Table 2}
\end{center}
\indent Throwing together these results we obtain
\[\{-\frac{n-1}{2}-\frac{1}{i} \mid i \in \mathbb{Z}_{>1}\} \subset \mathcal{P}_n.\]
\indent Now that we checked this expectation, we proceed proving
the theorem.
\\ \\
\emph{Remark.---} Notice that $m \in \mathcal{P}_{n-1}$ implies that $m \in
\mathcal{P}_{n}$. 
Indeed, any polynomial $f$ in $n-1$ variables can be considered as a
polynomial in $n$ variables. An embedded resolution for
$f^{-1}\{0\}\subset \mathbb{C}^{n-1}$ induces the obvious analogous
one for $f^{-1}\{0\}\subset \mathbb{C}^{n}=\mathbb{C}^{n-1}\times
\mathbb{C}$ and, since $\chi(\mathbb{C})=1$, the two associated
topological zeta functions are equal. From this observation it
follows that it is sufficient to prove that $[-(n-1)/2,-(n-2)/2)
\cap \mathbb{Q} \subset \mathcal{P}_n$.
As we showed in this section that $-(n-1)/2$ is contained in
$\mathcal{P}_{n-1}$ and thus in $\mathcal{P}_{n}$, we restrict
ourselves in the next sections to the subset $(-(n-1)/2,-(n-2)/2)
\cap \mathbb{Q}$.
\\ ${}$
\begin{center}
\textsc{3. The set $(-1/2,0)  \cap \mathbb{Q}$ is a subset of
$\mathcal{P}_2$}
\end{center} ${}$ \\
\indent Considering how candidate poles look like in the formula of
the topological zeta function written in terms of newton polyhedra
(see \cite{DenefLoeser1}), the number $-(b+2)/(2a+2b)$ seems to
appear as a candidate pole of the topological zeta function
associated to $f(x,y)=x^a(x^b+y^2)$, where $a$ and $b$ are
positive integers. An easy computation yields:\\
\\ \textbf{Lemma.---} \emph{When $a$ and $b$ run through $2
\mathbb{Z}_{>0}$, $a \not= 2$, the quotient $-(b+2)/(2a+2b)$ takes
all rational values in $(-1/2,0)$.}
\\ \\ \indent Taking the lemma into account, the functions $f(x,y)=x^a(x^b+y^2)$, where $a, b \in 2
\mathbb{Z}_{> 0}$ and $a \not= 2$, could be a pretty nice choice
to obtain all desired poles.
Easy calculations give the following dual resolution graph for $f$. \\ \\
\\
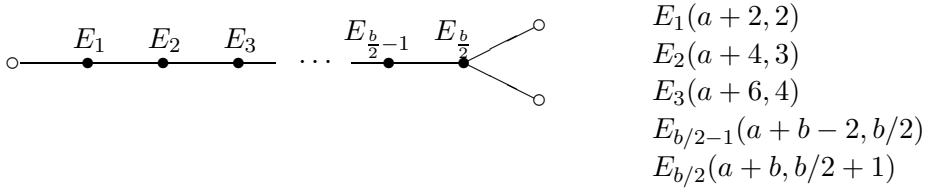
\begin{picture}(150,20)(-15,-5)
\put(-9,10){\line(1,0){8.5}} \put(-10,10){\circle{1.5}}
\put(0,10){\line(1,0){25}} \put(28,10){$\ldots$}
\put(35,10){\line(1,0){15}} \put(0,10){\circle*{1.5}}
\put(10,10){\circle*{1.5}} \put(20,10){\circle*{1.5}}
\put(40,10){\circle*{1.5}} \put(50,10){\circle*{1.5}}
\put(50,10){\line(2,1){9.3}} \put(50,10){\line(2,-1){9.3}}
\put(60,15){\circle{1.5}} \put(60,5){\circle{1.5}}
\put(-2,12){$E_1$} \put(8,12){$E_2$} \put(18,12){$E_3$}
\put(34,13){$E_{\frac{b}{2}-1}$} \put(46,13){$E_{\frac{b}{2}}$}
\put(75,15){$E_1(a+2,2)$} \put(75,10){$E_2(a+4,3)$}
\put(75,5){$E_3(a+6,4)$} \put(75,0){$E_{b/2-1}(a+b-2,b/2)$}
\put(75,-5){$E_{b/2}(a+b,b/2+1)$}
\end{picture}
\\ \\ Because $E_{b/2}$ is intersected three times by other components, Theorem $4.3$ in \cite{Veys2}
allows us to conclude that $-(b+2)/(2a+2b)$ is a pole of
$Z_{top,f}$.
\\  ${}$
\begin{center}
\textsc{4. The set $(-(n-1)/2,-(n-2)/2) \cap \mathbb{Q}$ is a subset of $\mathcal{P}_n,
\quad n \geq 3$}
\end{center} ${}$
\\ \indent As this set is a translation by $-1/2$ of expected
poles in dimension $n-1$, the Thom-Sebastiani principle in
\cite{DenefLoeser3} is again the motivation why we consider
\[f(x_1,\ldots,x_n)=x_n^2+ \cdots + x_3^2+ x_1^a(x_1^b+x_2^2),\]
where $a \in 2\mathbb{Z}_{> 0}$ and $a \not= 2$, to reach the set
$(-(n-1)/2,-(n-2)/2) \cap \mathbb{Q}$. ${}$ \\ \\
\emph{\textbf{Embedded resolution for \boldmath{$z^2+ x^a(x^b+y^2)$}}} \\
\indent Let us first explain in dimension $3$ which embedded
resolution we choose for $z^2+x^a(x^b+y^2)$ ($a,b \in
2\mathbb{Z}_{>0}$, $a \not=2$). We first blow up in the singular
locus $\{x=z=0\}$ of $f$ and further always in the singular locus
of the strict transform; the first $a/2$ times this is an affine
line and the last $b/2$ times it is a point. This is the special
case for $n=3$ in Table $3$.
\\The dual intersection graph looks as follows.\\ \\
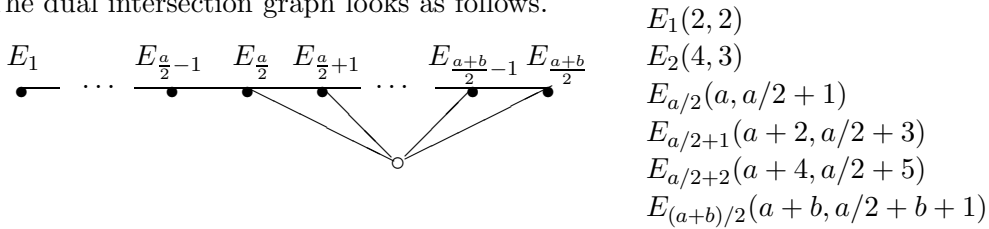
\begin{picture}(150,25)(-5,-5)
\put(0,15){\line(1,0){5}} \put(8,15){$\ldots$}
\put(15,15){\line(1,0){30}} \put(0,14.5){\circle*{1.5}}
\put(20,14.5){\circle*{1.5}} \put(30,14.5){\circle*{1.5}}
\put(40,14.5){\circle*{1.5}} \put(30,15){\line(2,-1){19.5}}
\put(40,15){\line(1,-1){9.15}} \put(50,5){\circle{1.5}}
\put(-2,18){$E_1$} \put(15,18){$E_{\frac{a}{2}-1}$}
\put(28,18){$E_{\frac{a}{2}}$} \put(36,18){$E_{\frac{a}{2}+1}$}
\put(47,15){$\ldots$} \put(51,6){\line(1,1){9}}
\put(60,14.5){\circle*{1.5}} \put(70,14.5){\circle*{1.5}}
\put(51,5.5){\line(2,1){19.5}} \put(55,15){\line(1,0){15}}
\put(54,18){$E_{\frac{a+b}{2}-1}$}
\put(67,18){$E_{\frac{a+b}{2}}$}
\put(83,23){$E_1(2,2)$} \put(83,18){$E_2(4,3)$}
\put(83,13){$E_{a/2}(a,a/2+1)$} \put(83,8){$E_{a/2+1}(a+2,a/2+3)$}
\put(83,3){$E_{a/2+2}(a+4,a/2+5)$}\put(83,-2){$E_{(a+b)/2}(a+b,a/2+b+1)$}
\end{picture}  \\
The candidate pole given by the last exceptional surface,
$E_{(a+b)/2}$, is equal to
\[-\frac{a/2+b+1}{a+b}=-\frac{b+2}{2a+2b}-\frac{1}{2},\]
and thus covers all rational numbers in $(-1,-1/2)$ if
$a$ and $b$ run over $2\mathbb{Z}_{>0}$ and $a \not= 2$.
\\ \\ \noindent \emph{\textbf{Embedded resolution for \boldmath{$x_n^2+ \cdots
+ x_3^2+ x_1^a(x_1^b+x_2^2)$}, $n > 3$}} \\
\indent The sequence of blow-ups in Table $3$ yields an embedded
resolution for \[f(x_1,\ldots,x_n)=x_n^2+ \cdots + x_3^2+
x_1^a(x_1^b+x_2^2),\] based on the previous one for $n=3$.
\\ \\
\begin{tabular}[t]{|c|c|c|}
\hline number $i$ of & center blow-up & equation strict transform
\\blow-up & & in relevant chart \\
\hline & & \\
$1$ & $x_1=x_3=x_4=\cdots=x_n=0$ & $x_n^2+ \cdots
+ x_3^2+ x_1^{a-2}(x_1^b+x_2^2)$  \\
$2$ & $x_1=x_3=x_4=\cdots=x_n=0$ & $x_n^2+ \cdots
+ x_3^2+ x_1^{a-4}(x_1^b+x_2^2)$  \\
$\vdots$ & $\vdots$ & $\vdots$ \\
$a/2$ & $x_1=x_3=x_4=\cdots=x_n=0$ & $x_n^2+ \cdots
+ x_3^2+ x_1^b+x_2^2$ \\
& & \\
$a/2+1$ & $(0,0,\ldots,0)$ & $x_n^2+ \cdots
+ x_3^2+ x_1^{b-2}+x_2^2$ \\
$a/2+2$ & $(0,0,\ldots,0)$ & $x_n^2+ \cdots
+ x_3^2+ x_1^{b-4}+x_2^2$  \\
$\vdots$ & $\vdots$ & $\vdots$ \\
$(a+b)/2$ & $(0,0,\ldots,0)$ & $x_n^2+ \cdots
+ x_3^2+ 1 +x_2^2$ \\
& & \\
\hline
\end{tabular}
\begin{center} \emph{Table 3} \end{center} ${}$\\
The dual intersection graph here looks as follows.
\\ 
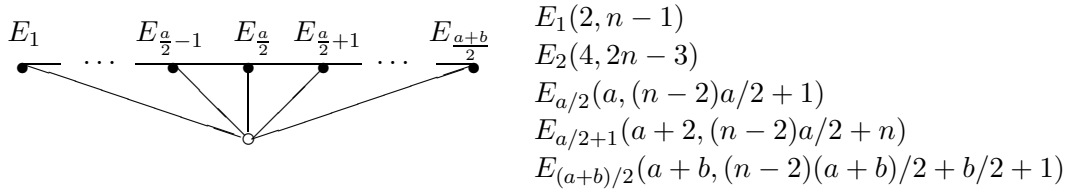
\begin{picture}(150,28)(-5,13)
\put(0,25){\line(3,-1){29}}
\put(20,25){\line(1,-1){9.5}} \put(30,25){\line(0,-1){9}}
\put(30.5,15){\line(1,1){9}} \put(30.5,15){\line(3,1){29}}
\put(30,15){\circle{1.5}} 
\put(0,25){\line(1,0){5}} \put(8,25){$\ldots$}
\put(15,25){\line(1,0){30}} \put(0,24.5){\circle*{1.5}}
\put(20,24.5){\circle*{1.5}} \put(30,24.5){\circle*{1.5}}
\put(40,24.5){\circle*{1.5}} \put(60,24.5){\circle*{1.5}}
\put(-2,28){$E_1$}
\put(15,28){$E_{\frac{a}{2}-1}$} \put(28,28){$E_{\frac{a}{2}}$}
\put(36,28){$E_{\frac{a}{2}+1}$} \put(47,25){$\ldots$}
\put(55,25){\line(1,0){5}}
\put(54,28){$E_{\frac{a+b}{2}}$} \put(68,30){$E_1(2,n-1)$}
\put(68,25){$E_2(4,2n-3)$} \put(68,20){$E_{a/2}(a,(n-2)a/2+1)$}
\put(68,15){$E_{a/2+1}(a+2,(n-2)a/2+n)$}
\put(68,10){$E_{(a+b)/2}(a+b,(n-2)(a+b)/2+b/2+1)$}
\end{picture}
\\ \\ Now $-\nu_{(a+b)/2}/N_{(a+b)/2}$ is equal to
\[-\frac{a/2+b+1 + ((a+b)/2)(n-3)}{a+b}=-\frac{b+2}{2a+2b}-\frac{n-2}{2},\]
which covers the interval $(-(n-1)/2,-(n-2)/2) \cap \mathbb{Q}$ when $a$ and $b$ vary in $2\mathbb{Z}_{>0}$ with $a \not=2$.
\\ \\ 
\noindent \emph{\textbf{The rational number \boldmath{$-\nu_{(a+b)/2}/N_{(a+b)/2}$} is a
pole of \boldmath{$Z_{top,f}$}}} \\
\indent For all $n \geq 3$ and $f(x_1,\ldots,x_n)=x_n^2+ \cdots +
x_3^2+ x_1^a(x_1^b+x_2^2)$, we calculate the residue of $Z_{top,f}$
in  $-\nu_{(a+b)/2}/N_{(a+b)/2}$. Observe that if $(a+b)/(2+b) \in
\mathbb{Z}$, the exceptional variety $E_{(a+b)/(2+b)}$ induces the
same candidate pole as $E_{(a+b)/2}$. The other exceptional
varieties always give rise to
other candidate poles. \\
The subsets playing a role in the contribution of
$E_{(a+b)/(2+b)}$ to the residue are $J_1=\{\frac{a+b}{2+b}\},
\quad J_2=\{\frac{a+b}{2+b},\frac{a+b}{2+b}-1\},
J_3=\{\frac{a+b}{2+b},
 \frac{a+b}{2+b}+1\}, J_4=\{\frac{a+b}{2+b},0\}, \quad J_5=\{\frac{a+b}{2+b},\frac{a+b}{2+b}-1,0\},
 \quad J_6=\{\frac{a+b}{2+b},\frac{a+b}{2+b}+1,0\}.$
Notice that when $n=3$, $E_{(a+b)/(2+b)}$ does not intersect
$E_0$.\\ \indent We have that $E_{(a+b)/(2+b)}$ is isomorphic to the
cartesian product of $\mathbb{A}^1$ and the blowing-up of
$\mathbb{P}^{n-2}$ in a point. It is also easy to describe the whole
intersection configuration on $E_{(a+b)/(2+b)}$.
\begin{center}
\begin{tabular}[t]{|c|c|c|} \hline $\chi(E^\circ_{J_k})$ & n
odd & n even
\\ \hline & & \\
$k=1$ & $0$ & $0$  \\
$k=2$ & $1$ & $0$\\
$k=3$ & $1$ & $0$ \\
$k=4$ & $0$ & $0$ \\
$k=5$ & $n-3$ & $n-2$ \\
$k=6$ & $n-3$ & $n-2$ \\
& & \\
\hline
\end{tabular}
\end{center}
\begin{center}
\emph{Table 4}
\end{center}
With the relevant Euler--Poincar\'e characteristics of Table $4$
and $\alpha_{(a+b)/(2+b)-1}=1/i$, $\alpha_{(a+b)/(2+b)+1}=-1/i$,
we see that $E_{(a+b)/(2+b)}$ does not give any contribution to
the residue in $-\nu_{(a+b)/2}/N_{(a+b)/2}$. Alternatively, this
is implied by \cite[Proposition 6.5]{Veysconfigurations}. This
means we only have to
take the contribution of $E_{(a+b)/2}$ into account.\\
\indent To compute this contribution the relevant subsets for the
summation in the formula of the topological zeta function are
$I_1=\{\frac{a+b}{2}\}, \mbox{ }
 I_2=\{\frac{a+b}{2},\frac{a+b}{2}-1\}, \mbox{ }
 I_3=\{\frac{a+b}{2},0\}, \mbox{ }
 I_4=\{\frac{a+b}{2},\frac{a+b}{2}-1,0\}.$
The Euler--Poincar\'e characteristics $\chi(E^\circ_{I_j})$, $1
\leq j \leq 4$, are the same as those given in Table $1$ and we have
$\alpha_0 =-((n-4)a+(n-3)b+2)/(2(a+b))$ and $\alpha_{(a+b)/2-1}=(2-a)/(a+b)$.
\\ \indent As the residue then is equal to
\[\frac{(-2+3a+2b)(na-2a-b+nb+2)}{(-2+a)(a+b)(na-4a+2+nb-3b)}
 \quad \mbox{ for $n$ odd and}\]
\[\frac{(2 + b)(n a - 2 a - b + n b + 2)}{(-2 + a)(a + b)(n a - 4 a + 2 + n b - 3 b)} \quad \mbox{ for $n$ even},\]
we find that $-(\nu_{(a+b)/2})/(N_{(a+b)/2})=-(b+2)/(2a+2b)-(n-2)/2$
is a pole of $Z_{top,f}$. \\ \\ \indent We conclude that
$(-(n-1)/2,-(n-2)/2)  \cap \mathbb{Q}\subset \mathcal{P}_n$, for all $n \geq 3$.
\\ ${}$
 \begin{center}
\textsc{5. Some remarks}
\end{center} ${}$ \\
\textbf{($1$)} Instead of achieving this result with the method of
resolution of singularities one can find the poles of the
topological zeta function of the polynomials
\[x_n^2+ \cdots
+ x_3^2+ x_1^a(x_1^b+x_2^2) \quad \mbox{ and } \quad x_n^2+ \cdots +
x_2^2+ x_1^i \] with the help of Newton polyhedra. Indeed, we can
write down the topolo\-gical zeta function for these polynomials
using the formula of Denef and Loeser in \cite{DenefLoeser1}. For
example if $f(x_1,\ldots,x_n)=x_n^2+ \cdots + x_3^2+
x_1^a(x_1^b+x_2^2)$, where $a$ and $b$ are positive even
integers and $a\not=2$, put $A:=(a+b)s+1+b/2+(n-2)(a+b)/2$ and
$B:=as+1+(n-2)a/2$. We get
\begin{eqnarray*}
Z_{top,f}(s) & = & (n-1)\frac{b}{2AB}+\frac{1}{A}+(n-2)\frac{a}{2B}\\
& & + \frac{s}{s+1}\left(\sum_{d=1}^{n-1}{  n-2 \choose
  d+1}
\left(\frac{a}{2B}+\frac{b}{2AB}\right)(-2)^d \right. \\
& & +\left. \sum_{d=1}^{n-1}{  n-1 \choose d}
\frac{1}{A}(-2)^d+\sum_{d=1}^{n-2}{  n-2 \choose   d}
\frac{b}{2AB}(-2)^d\right).
\end{eqnarray*}
Handling the problem in this way leads to the same results. One
just has to be careful with the dual cones of some faces, namely
those that are not a rational simplicial cone.
\\ \\ 
\textbf{($2$)} With a similar definition of $\mathcal{P}_n$ in each
case, the same results hold
 for local and global versions of the motivic zeta function, the Hodge zeta function and
 Igusa's zeta function. Indeed, the results for the topological zeta function imply the results
 for those `finer' zeta functions.

\footnotesize{

 \end{document}